\documentclass[a4paper,12pt]{elsart}
\usepackage{amsmath,amssymb,latexsym}
\newcommand{\fd}{\mathbb{F}}

\newcommand{\Z}{\mathbb{Z}}

\newtheorem{theorem}{Theorem}
\newtheorem{lemma}[theorem]{Lemma}
\newtheorem{corollary}[theorem]{Corollary}
\newtheorem{proposition}[theorem]{Proposition}
\newcommand{\dis}{\displaystyle}

\begin{document}

\begin{frontmatter}

\title{Infinite Families of Recursive Formulas Generating Power Moments of Kloosterman Sums: $O^-(2n,2^r)$ Case}

\author{Dae San Kim}

\address{Department of Mathematics, Sogang University, Seoul 121-742, South Korea}

\begin{abstract}
In this paper, we construct eight infinite families of binary linear
codes  associated with double cosets with respect to certain
maximal parabolic subgroup  of the special orthogonal group
$SO^-(2n,2^r)$. Then we obtain four infinite families of recursive
formulas for the power moments of Kloosterman sums
  and  four those of 2-dimensional Kloosterman sums in terms of the frequencies of weights in the codes.
  This is done via Pless power moment identity and by utilizing the explicit expressions of exponential sums over those double
  cosets  related to the evaluations of "Gauss sums" for the orthogonal  groups $O^-(2n,2^r)$
\end{abstract}

\thanks{Email address : dskim@sogang.ac.kr.}


\end{frontmatter}

\section{Introduction}

Let  $\psi$ be a nontrivial additive character of the finite field
$\fd_q$  with $q = p^r$ elements ($p$ a prime), and let  $m$ be a
positive integer. Then the $m$-dimensional  Kloosterman sum
$K_m(\psi;a)$ (\cite{RH}) is defined by
\[
K_m(\psi;a) = \sum_{\alpha_1, \ldots, \alpha_m \in \fd_q^*}
\psi(\alpha_1 + \cdots + \alpha_m +
a\alpha_1^{-1}\cdots\alpha_m^{-1}) (a \in \fd_q^*).
\]
In particular, if $m=1$, then  $K_1(\psi;a)$ is simply denoted by
$K(\psi;a)$, and is called the Kloosterman sum. The Kloosterman sum
was introduced in 1926 to give an estimate for the Fourier
coefficients of modular forms (cf. \cite{HD}, \cite{JH}). It has
also been studied to solve  various problems in coding theory and
cryptography over finite fields of characteristic two (cf.
\cite{PTV}, \cite{HPT}).

For each nonnegative integer $h$, by $MK_m(\psi)^h$ we will denote
the $h$-th moment of the $m$-dimensional Kloosterman sum
$K_m(\psi;a)$. Namely, it is given by
\[
MK_m(\psi)^h = \sum_{a \in \fd_q^*} K_m(\psi;a)^h.
\]
If  $\psi = \lambda$ is the canonical additive character of $\fd_q$,
then $MK_m(\lambda)^h$ will be simply denoted by $MK_m^h$. If
further $m=1$, for brevity $MK_1^h$ will be indicated by $MK^h$.

Explicit computations on power moments of Kloosterman sums were
begun with the paper \cite{HS} of Sali\'{e} in 1931, where he
showed, for any odd prime $q$,
\[
MK^h = q^2M_{h-1} - (q-1)^{h-1} + 2(-1)^{h-1} ~~~(h \geq 1).
\]
Here $M_0=0$, and for $h \in \Z_{>0}$,
\[
M_h = \Big | \Big \{(\alpha_1,\ldots,\alpha_h) \in (\fd_q^*)^h
\mid \sum_{j=1}^h \alpha_j = 1 = \sum_{j=1}^h \alpha_j^{-1} \Big
\} \Big |.
\]
For  $q = p$ odd prime, Sali\'{e} obtained  $MK^1, MK^2, MK^3, MK^4$
in \cite{HS} by determining $M_1, M_2, M_3$. $MK^{5}$ can be
expressed in terms of the $p$-th eigenvalue for a weight 3 newform
on $\Gamma_{0}$(15) (cf. \cite{RL}, \cite{CJM}). $MK^{6}$ can be
expressed in terms of the $p$-th eigenvalue for a weight 4 newform
on $\Gamma_{0}$(6) (cf. \cite{KJBD}). Also, based on numerical
evidence, in \cite{RJE} Evans was led to propose a conjecture which
expresses $MK^7$ in terms of Hecke eigenvalues for a weight 3
newform on $\Gamma_0 (525)$ with quartic nebentypus of conductor
105. For more details about this brief history of explicit
computations on power moments of Kloosterman sums, one is referred
to Section IV of \cite{Kim1}.

From now on, let us assume that $q = 2^r$. Carlitz\cite{L1}
evaluated $MK^h$ for the other values of $h$ with $h \leq 10$
(cf.\cite{M1}). Recently, Moisio was able to find explicit
expressions of $MK^h$, for $h \leq 10$ (cf. \cite{M1}). This was
done, via Pless power moment identity, by connecting moments of
Kloosterman sums and the frequencies of weights in the binary
Zetterberg code of length $q+1$, which were known by the work of
Schoof and Vlugt in \cite{RS}.

In \cite{Kim1}, the binary linear codes $C(SL(n,q))$ associated with
finite special linear groups $SL(n,q)$  were constructed when $n,q$
are both powers of two. Then obtained was a recursive formula for
the power moments of multi-dimensional Kloosterman sums in terms of
the frequencies of weights in $C(SL(n,q))$. In particular, when
$n=2$, this gives a recursive formula for the power moments of
Kloosterman sums. Also, in order to get recursive formulas for the
power moments of Kloosterman and 2-dimensional Kloosterman sums, we
constructed in \cite{Kim2} three binary linear codes $C(SO^+(2,q))$,
$C(O^+(2,q))$, $C(SO^+(4,q))$, respectively associated with
$SO^+(2,q)$, $O^+(2,q)$, $SO^+(4,q)$, and in \cite{Kim3} three
binary linear codes $C(SO^-(2,q))$, $C(O^-(2,q))$, $C(SO^-(4,q))$,
respectively associated with $SO^-(2,q)$, $O^-(2,q)$, $SO^-(4,q)$.
All of these were done via Pless power moment identity and by
utilizing our previous results on explicit expressions of Gauss sums
for the stated finite classical groups. So, all in all, we had only
a handful of recursive formulas generating power moments of
Kloosterman and 2-dimesional Kloosterman sums.

In this paper, we will be able to produce four infinite families
of recursive formulas generating power moments of Kloosterman sums
and four those of  2-dimensional Kloosterman sums. To do that, we
construct eight  infinite families of binary linear codes
$C(DC_1^+(n,q))$ $(n=2,4,\ldots)$, $C(DC_1^-(n,q))$
$(n=1,3,\ldots)$, both associated with $Q^-\sigma_{n-1}^-Q^-$ ;
$C(DC_2^+(n,q))$ $(n=2,4,\ldots)$, $C(DC_2^-(n,q))$
$(n=3,5,\ldots)$, both associated with $Q^-\sigma_{n-2}^-Q^-$;
$C(DC_3^+(n,q))$ $(n=2,4,\ldots)$, $C(DC_3^-(n,q))$
$(n=3,5,\ldots)$, both associated with $\rho
Q^-\sigma_{n-2}^-Q^-$; $C(DC_4^+(n,q))$ $(n=4,6,\ldots)$,
$C(DC_4^-(n,q))$ $(n=3,5,\ldots)$, both associated with $\rho
Q^-\sigma_{n-3}^-Q^-$, with respect to the maximal parabolic
subgroup  $Q^- = Q^-(2n,q)$ of the special orthogonal group
$SO^-(2n,q)$, and express those power moments in terms of the
frequencies of weights in each code. Then, thanks to our previous
results on the explicit expressions of exponential sums over those
double cosets related to the evaluations of ``Gauss sums" for the
orthogonal groups $O^-(2n,q)$\cite{Kim8}, we can express the
weight of each codeword in the duals of the codes in terms of
Kloosterman or 2-dimensional Kloosterman sums. Then our formulas
will follow immediately from the Pless power moment identity.
Analogously to these, in \cite{Kim4} and \cite{Kim5}, we obtained
infinite families of recursive formulas for power moments of
Kloosterman and 2-dimensional Kloosterman sums by constructing
binary codes associated with double cosets with respect to certain
maximal parabolic subgroup of the symplectic group $Sp(2n,q)$ and
the orthogonal group $O^+(2n,q)$, respectively.

Theorem 1 in the following(cf. (17), (18), (20)-(25)) is the main
result of this paper. Henceforth, we agree that the binomial
coefficient $\binom{b}{a}=0$, if $a > b$ or $a<0$. To simplify
notations, we introduce the following ones which will be used
throughout this paper at various places.

\begin{align}
A_1^+(n,q) &= q^{\frac{1}{4}(5n^2-2n-4)} (q^{n-1} - 1)
\prod_{j=1}^{\frac{(n-2)}{2}} (q^{2j-1} -1),\\
B_1^+(n,q) &= (q+1)q^{\frac{1}{4}n^2}
\prod_{j=1}^{\frac{(n-2)}{2}} (q^{2j}-1),\\
A_2^+(n,q) &= q^{\frac{1}{4}(5n^2-2n-8)} \left[ \substack{n-1\\1}
\right]_q \prod_{j=1}^{\frac{(n-2)}{2}} (q^{2j-1} -1),
\end{align}
\begin{align}
B_2^+(n,q) &=
(q+1)q^{\frac{1}{4}(n-2)^2}(q^{n-1}-1)\prod_{j=1}^{\frac{(n-2)}{2}}
(q^{2j}-1),\\
A_3^+(n,q) &= (q+1)q^{\frac{1}{4}(5n^2-2n-8)} \left[
\substack{n-1\\1} \right]_q \prod_{j=1}^{\frac{(n-2)}{2}}
(q^{2j-1} -1),\\
B_3^+(n,q) &=
q^{\frac{1}{4}(n-2)^2}(q^{n-1}-1)\prod_{j=1}^{\frac{(n-2)}{2}}
(q^{2j}-1),\\
A_4^+(n,q) &= (q+1)q^{\frac{1}{4}(5n^2-6n-4)} \left[
\substack{n-1\\2} \right]_q \prod_{j=1}^{\frac{(n-2)}{2}}
(q^{2j-1} -1),\\
B_4^+(n,q) &=
q^{\frac{1}{4}(n-2)^2}(q^{n-1}-1)\prod_{j=1}^{\frac{(n-2)}{2}}
(q^{2j}-1),\\
A_1^-(n,q) &= q^{\frac{5}{4}(n^2-1)} \prod_{j=1}^{\frac{(n-1)}{2}} (q^{2j-1} -1),\\
B_1^-(n,q) &= (q+1)q^{\frac{1}{4}(n-1)^2}
\prod_{j=1}^{\frac{(n-1)}{2}} (q^{2j}-1),\\
A_2^-(n,q) &= q^{\frac{1}{4}(5n^2-4n-5)} \left[ \substack{n-1\\1}
\right]_q \prod_{j=1}^{\frac{(n-1)}{2}}(q^{2j-1} -1),\\
B_2^-(n,q) &=
(q+1)q^{\frac{1}{4}(n-1)^2}\prod_{j=1}^{\frac{(n-1)}{2}}
(q^{2j}-1),\\
A_3^-(n,q) &= (q+1)q^{\frac{1}{4}(5n^2-4n-5)} \left[
\substack{n-1\\1}
\right]_q \prod_{j=1}^{\frac{(n-1)}{2}}(q^{2j-1} -1),\\
B_3^-(n,q) &= q^{\frac{1}{4}(n-1)^2}\prod_{j=1}^{\frac{(n-1)}{2}}
(q^{2j}-1),\\
A_4^-(n,q) &= (q+1)q^{\frac{1}{4}(5n^2-4n-9)} \left[
\substack{n-1\\2}
\right]_q \prod_{j=1}^{\frac{(n-3)}{2}}(q^{2j-1} -1),\\
B_4^-(n,q) &= q^{\frac{1}{4}(n-3)^2}(q^{n-2}-1)(q^{n-1}-1)
\prod_{j=1}^{\frac{(n-3)}{2}} (q^{2j}-1).
\end{align}
\emph{From now on, it is assumed that either +signs or - signs are
chosen everywhere, whenever $\pm$ signs appear. }

\begin{theorem}
Let $q = 2^r$. Then, with the notations in (1)-(16), we have the
following. \\
(a) With $i=1$  and $+$ signs everywhere for $\pm$ signs,  we have
a recursive formula generating power moments of Kloosterman sums
over $\fd_q$, for each $n \geq 2$ even and all $q$; with $i=3$ and
$+$ signs everywhere for $\pm$ signs,  we have such a formula, for
either each $n \geq 4$ even and all $q$, or $n=2$ and $q \geq 8$;
with $i=1$ and $-$ signs everywhere for $\pm$ signs, we have such
a formula, for each $n \geq 1$ odd and all $q$;  with $i=3$ and
$-$ signs everywhere for $\pm$ signs,  we have such a formula, for
each $n \geq 3$ odd and all $q$.
\begin{multline}
(\pm(-1))^h MK^h = - \sum_{l=0}^{h-1} (\pm(-1))^l {h \choose l}
B_i^{\pm}(n,q)^{h-l}MK^l + qA_i^{\pm}(n,q)^{-h} \\
 \times \sum_{j=0}^{min\{N_i^{\pm}(n,q),h\}} (-1)^j
C_{i,j}^{\pm}(n,q) \sum_{t=j}^h t! S(h,t) 2^{h-t}
{N_i^{\pm}(n,q)-j \choose N_i^{\pm}(n,q)-t} (h=1,2,\ldots),
\end{multline}
where $N_i^{\pm}(n,q) = |DC_i^{\pm}(n,q)| =
A_i^{\pm}(n,q)B_i^{\pm}(n,q)$, and
$\{C_{i,j}^{\pm}(n,q)\}_{j=0}^{N_i^{\pm}(n,q)}$ is the weight
distribution of the binary code $C(DC_i^{\pm}(n,q))$ given by
\begin{multline}
C_{i,j}^{\pm}(n,q) = \sum {q^{-1}A_i^{\pm}(n,q)(B_i^{\pm}(n,q)\pm
1) \choose \nu_0} \\
\times \prod_{tr(\beta^{-1})=0}
{q^{-1}A_i^{\pm}(n,q)(B_i^{\pm}(n,q)\pm (q+1)) \choose \nu_\beta}
\prod_{tr(\beta^{-1})=1} {q^{-1}A_i^{\pm}(n,q)(B_i^{\pm}(n,q)\pm
(-q+1)) \choose \nu_\beta},
\end{multline}
where the sum is over all the sets of nonnegative integers
$\{\nu_\beta\}_{\beta \in \fd_q}$ satisfying $\dis\sum_{\beta \in
\fd_q} \nu_\beta = j$ and $\dis\sum_{\beta \in \fd_q} \nu_\beta
\beta = 0$. In addition, $S(h,t)$ is the Stirling number of the
second kind defined by
\begin{equation}
S(h,t) = \frac{1}{t!}\sum_{j=0}^t (-1)^{t-j} {t \choose j} j^h.
\end{equation}\\
(b) With $+$ signs everywhere for $\pm$ signs,  we have recursive
formulas generating power moments of 2-dimensional Kloosterman
sums over $\fd_q$ and even power moments of Kloosterman sums over
$\fd_q$, for each $n \geq 2$ even and  $q \geq 4$; with $-$ signs
everywhere for $\pm$ signs, we have such formulas, for each $n
\geq 3$ odd and $q \geq 4$.
\begin{multline}
(\pm1)^h MK_2^h = - \sum_{l=0}^{h-1} (\pm1)^l {h \choose l}
(B_2^{\pm}(n,q)\pm q)^{h-l}MK_2^l + qA_2^{\pm}(n,q)^{-h} \\
 \times \sum_{j=0}^{min\{N_2^{\pm}(n,q),h\}} (-1)^j
C_{2,j}^{\pm}(n,q) \sum_{t=j}^h t! S(h,t) 2^{h-t}
{N_2^{\pm}(n,q)-j \choose N_2^{\pm}(n,q)-t} (h=1,2,\ldots),
\end{multline}
and
\begin{multline}
(\pm1)^h MK^{2h} = - \sum_{l=0}^{h-1} (\pm1)^l {h \choose l}
B_2^{\pm}(n,q)^{h-l}MK^{2l} + qA_2^{\pm}(n,q)^{-h} \\
 \times \sum_{j=0}^{min\{N_2^{\pm}(n,q),h\}} (-1)^j
C_{2,j}^{\pm}(n,q) \sum_{t=j}^h t! S(h,t) 2^{h-t}
{N_2^{\pm}(n,q)-j \choose N_2^{\pm}(n,q)-t} (h=1,2,\ldots),
\end{multline}
where $N_2^{\pm}(n,q) = |DC_2^{\pm}(n,q)| =
A_2^{\pm}(n,q)B_2^{\pm}(n,q)$, and
$\{C_{2,j}^{\pm}(n,q)\}_{j=0}^{N_2^{\pm}(n,q)}$ is the weight
distribution of the binary code $C(DC_2^{\pm}(n,q))$ given by
\begin{multline}
C_{2,j}^{\pm}(n,q) = \sum {q^{-1}A_2^{\pm}(n,q)(B_2^{\pm}(n,q)\pm
(q+1-q^2)) \choose \nu_0} \\
\times \prod_{\substack{|\tau|<2\sqrt{q}\\\tau \equiv -1(mod4)}}
\prod_{K(\lambda;\beta^{-1})=\tau}
{q^{-1}A_2^{\pm}(n,q)(B_2^{\pm}(n,q)\pm (q+1-q\tau)) \choose
\nu_\beta},
\end{multline}
where the sum is over all the sets of nonnegative integers
$\{\nu_\beta\}_{\beta \in \fd_q}$ satisfying $\dis\sum_{\beta \in
\fd_q} \nu_\beta = j$ and $\dis\sum_{\beta \in \fd_q} \nu_\beta
\beta = 0$.\\

(c)  With $+$ signs everywhere for $\pm$ signs,  we have recursive
formulas generating power moments of 2-dimensional Kloosterman
sums over $\fd_q$ and even power moments of Kloosterman sums over
$\fd_q$, for each $n \geq 4$ even and $q \geq 4$;  with $-$ signs
everywhere for $\pm$ signs, we have such formulas, for each $n
\geq 3$ odd and $q \geq 4$.
\begin{multline}
(\pm 1)^h MK_2^h = - \sum_{l=0}^{h-1} (\pm 1)^l {h \choose l}
\{B_4^{\pm}(n,q) \pm q^2 \}^{h-l} MK_2^l + qA_4^{\pm}(n,q)^{-h} \\
\times \sum_{j=0}^{min \{N_4^{\pm}(n,q),h\}} (-1)^j
C_{4,j}^{\pm}(n,q) \sum_{t=j}^h t!S(h,t) 2^{h-t} {N_4^{\pm}(n,q)-j
\choose N_4^{\pm}(n,q)-t} (h=1,2,\ldots),
\end{multline}
and
\begin{multline}
(\pm 1)^h MK^{2h} = - \sum_{l=0}^{h-1} (\pm 1)^l {h \choose l}
\{B_4^{\pm}(n,q) \pm (q^2-q) \}^{h-l} MK^{2l} + qA_4^{\pm}(n,q)^{-h} \\
\times \sum_{j=0}^{min \{N_4^{\pm}(n,q),h\}} (-1)^j
C_{4,j}^{\pm}(n,q) \sum_{t=j}^h t!S(h,t) 2^{h-t} {N_4^{\pm}(n,q)-j
\choose N_4^{\pm}(n,q)-t} (h=1,2,\ldots),
\end{multline}
where $N_4^{\pm}(n,q) = |DC_4^{\pm}(n,q)| = A_4^{\pm}(n,q)
B_4^{\pm}(n,q)$, and
$\{C_{4,j}^{\pm}(n,q)\}_{j=0}^{N_4^{\pm}(n,q)}$ is the weight
distribution of the binary code $C(DC_4^{\pm}(n,q))$ given by
\begin{multline}
C_{4,j}^{\pm}(n,q) = \sum {q^{-1} A_4^{\pm}(n,q)(B_4^{\pm}(n,q)
\pm (q^2+1-q^3)) \choose \nu_0}\\
\times \prod_{\substack{|\tau|<2\sqrt{q}\\\tau \equiv -1(mod 4)}}
\prod_{K(\lambda;\beta^{-1})=\tau} {q^{-1}
A_4^{\pm}(n,q)(B_4^{\pm}(n,q) \pm (q^2+1-q\tau)) \choose
\nu_\beta},
\end{multline}
where the sum is over all the sets of nonnegative integers
$\{\nu_\beta\}_{\beta \in \fd_q}$ satisfying $\dis\sum_{\beta \in
\fd_q} \nu_\beta = j$, and $\dis \sum _{\beta \in \fd_q} \nu_\beta
\beta = 0$.
\end{theorem}

The following corollary is just the  $n=2$ and  $n=1$ cases of (a)
in the above. It is amusing to note that the recursive formula in
(26) and (27), obtained from the binary code $C(DC_1^-(1,q))$
associated with the double coset $DC_1^-(1,q) = Q^-(2,q)$, is the
same as the one in (\cite{Kim3}, (1), (2)), gotten from the binary
code $C(SO^-(2,q))$ associated with the special orthogonal group
$SO^-(2,q)$.

\begin{corollary}
(a) For all $q$, and $h=1,2,\ldots$,
\begin{multline*}
MK^h = \sum_{l=0}^{h-1} (-1)^{h+l+1} {h \choose l} (q^2+q)^{h-l}
MK^l + q^{1-3h} (q-1)^{-h}\\
\times \sum_{j=0}^{min \{q^4(q^2-1),h\}} (-1)^{h+j}
C_{1,j}^{+}(2,q) \sum_{t=j}^h t! S(h,t) 2^{h-t} {q^4(q^2-1)-j
\choose q^4(q^2-1)-t},
\end{multline*}
where $\{C_{1,j}^+(2,q)\}_{j=0}^{q^4(q^2-1)}$ is the weight
distribution of $C(DC_1^+(2,q))$ given by
\begin{multline*}
C_{1,j}^+(2,q) = \sum {q^2(q-1)(q^2+q+1) \choose \nu_0}\\
\times \prod_{tr(\beta^{-1})=0} {q^2(q+1)(q^2-1) \choose \nu_\beta
} \prod_{tr(\beta^{-1})=1} {q^2(q-1)(q^2+1) \choose \nu_\beta}.
\end{multline*}
Here the sum is over all the sets of nonnegative integers
$\{\nu_\beta\}_{\beta \in \fd_q}$ satisfying $\dis\sum_{\beta \in
\fd_q} \nu_\beta = j$, and $\dis \sum _{\beta \in \fd_q} \nu_\beta
\beta = 0$. In addition, $S(h,t)$ is the Stirling number of the
second kind as defined in (19).\\

(b) For all $q$, and $h=1,2,\ldots$,
\begin{multline}
MK^h = - \sum_{l=0}^{h-1} {h \choose l} (q+1)^{h-l} MK^l\\
+ q \sum_{j=0}^{min \{q+1, h\}} (-1)^j C_{1,j}^-(1,q) \sum_{t=j}^h
t! S(h,t) 2^{h-t} {q+1-j \choose q+1-t},
\end{multline}
where $\{C_{1,j}^-(1,q)\}_{j=0}^{q+1}$ is the weight distribution
of $C(DC_1^-(n,q))$ given by
\begin{equation}
C_{1,j}^-(n,q) = \sum {1 \choose \nu_0} \prod_{tr(\beta^{-1})=1}
{2\choose \nu_\beta}.
\end{equation}
Here the sum is  over all the sets of nonnegative integers
$\{\nu_0\} \cup \{\nu_\beta\}_{tr(\beta^{-1})=1}$ satisfying
$\nu_0 + \dis\sum _{tr(\beta^{-1})=1} \nu_\beta = j$ and $\dis\sum
_{tr(\beta^{-1})=1} \nu_\beta \beta = 0$.
\end{corollary}

\section{$O^-(2n,q)$}

 For more details about the results of this section, one is referred to the paper
\cite{Kim8}. Throughout this paper, the following notations will
be used:
\begin{itemize}
 \item [] $q = 2^r$ ($r \in \Z_{>0}$),\\
 \item [] $\fd_q$ = the finite field with $q$ elements,\\
 \item [] $Tr A$ = the trace of $A$ for a square matrix $A$,\\
 \item [] $^tB$ = the transpose of $B$ for any matrix $B$.
\end{itemize}\

 Let $\theta^-$ be the nondegenerate quadratic form on the vector
 space $\fd_q^{2n \times 1}$ of all $2n \times 1$ column vectors over $\fd_q$, given by
 \begin{equation}
\theta^-(\sum_{i=1}^{2n} x_i e^i) = \sum_{i=1}^{n-1} x_i x_{n-1+i}
+ x_{2n-1}^2 + x_{2n-1}x_{2n} + ax_{2n}^2,
\end{equation}
  where \{$e^1=^t[10\ldots0], e^2=^t[01\ldots0],\ldots,e^{2n}=^t[0\ldots01]$\} is the standard basis of $\fd_q^{2n \times
  1}$, and $a$ is a fixed element in $\fd_q$ such that $z^2+z+a$
  is irreducible over $\fd_q$, or equvalently $a \in \fd_q
  \backslash \Theta (\fd_q)$, where $\Theta(\fd_q) = \{\alpha^2 + \alpha \mid \alpha \in \fd_q
  \}$ is a subgroup of index $2$ in the additive group $\fd_q^+$
  of $\fd_q$.

Let $\delta_a$(with $a$ in the above paragraph), $\eta$ denote
respectively the  $2 \times 2$ matrices over $\fd_q$ given by:
\[
\delta_a = \left[%
\begin{array}{cc}
  1 & 1 \\
  0 & a \\
\end{array}%
\right], \eta = \left[%
\begin{array}{cc}
  0 & 1 \\
  1 & 0 \\
\end{array}%
\right].
\]
Then the group $O^-(2n,q)$ of all isometries of $(\fd_q^{2n \times
1},\theta^-)$ consists of all matrices
\[
\left[%
\begin{array}{ccc}
  A & B & e \\
  C & D & f \\
  g & h & i \\
\end{array}%
\right] (A,B,C,D ~~(n-1)\times(n-1), e,f ~~(n-1) \times 2, g,h ~~2
\times (n-1))
\] in $GL(2n,q)$ satisfying the relations:
\begin{align*}
            &{}^tAC + {}^tg\delta_a g ~~\text{is alternating,}\\
            &{}^tBD + {}^th\delta_a h ~~\text{is alternating,}\\
            &{}^tef + {}^ti \delta_a i + \delta_a ~~\text{is
            alternating,}\\
            &{}^tAD + {}^tCB + {}^tg\eta h = 1_{n-1},\\
            &{}^tAf + {}^tCe + {}^tg\eta i = 0,\\
            &{}^tBf + {}^tDe + {}^th\eta i = 0.
\end{align*}
Here an  $n \times n$ matrix $(a_{ij})$ is called alternating if
\[
\left\{%
\begin{array}{ll}
    a_{ii} = 0, & \hbox{for $1 \leq i \leq n$,} \\
    a_{ij} = -a_{ji} = a_{ji}, & \hbox{for $1 \leq i < j \leq n$.} \\
\end{array}%
\right.
\]

$P^- = P^-(2n,q)$  is the maximal parabolic subgroup of
$O^-(2n,q)$ defined by:
\[
P^-(2n,q) = \left\{ \left[%
\begin{array}{ccc}
  A & 0 & 0\\
  0 & ^tA^{-1} & 0 \\
  0 & 0 & i \\
\end{array}%
\right] \left[%
\begin{array}{ccc}
  1_{n-1} & B & ^th^ti\eta i \\
  0 & 1_{n-1} & 0 \\
  0 & h & 1_2 \\
\end{array}%
\right] \Biggl \lvert %
\begin{array}{c}
  A \in GL(n-1,q), i \in O^{-}(2,q),\\
  {}^tB + {}^th \delta_a h ~~\text{is alternating}\\
\end{array}%
\right\},
\]
where $O^-(2,q)$ is the group of all isometries of $(\fd_q^{2 \times 1}, \theta^-)$ with
\[
\theta^-(x_1e^1 + x_2e^2) = x_1^2 + x_1 x_2 + ax_2^2 ~(cf. (28)).
\]
One can show that
\begin{equation}
O^-(2,q) = SO^-(2,q) \coprod \left[%
\begin{array}{cc}
  1 & 1 \\
  0 & 1 \\
\end{array}%
\right] SO^-(2,q),
\end{equation}
\begin{align*}
SO^-(2,q) &= \left\{ \left[%
\begin{array}{cc}
  d_1 & ad_2 \\
  d_2 & d_1+d_2 \\
\end{array}%
\right] \Big| d_1^2 + d_1 d_2 + ad_2^2 = 1
             \right\}\\
             &= \left\{ \left[%
\begin{array}{cc}
  d_1 & ad_2 \\
  d_2 & d_1+d_2 \\
\end{array}%
\right] \Big| %
\begin{array}{c}
 d_1 + d_2 b \in \fd_q(b), with\\
 N_{ \fd_q(b)/\fd_q}(d_1+d_2b) = 1
 \end{array}%
             \right\},
\end{align*}
where  $b \in \overline{\fd_q}$ is a root of the irreducible
polynomial $z^2+z+a$ over $\fd_q$.  $SO^-(2,q)$ is a subgroup of
index 2 in $O^-(2,q)$, and
\[
|SO^-(2,q)| = q+1, |O^-(2,q)| = 2(q+1).
\]
 $SO^-(2,q)$ here is defined as the kernel of a certain epimorphism $\delta^{-} : O^-(2n,q) \rightarrow \fd_2^+$ (cf. \cite{Kim8}, (3.45)).

The Bruhat decomposition of $O^-(2n,q)$  with respect to $P^- =
P^-(2n,q)$  is
\begin{equation}
O^-(2n,q) = \coprod_{r=0}^{n-1} P^- \sigma_r^- P^-,
\end{equation}
where
\[
\sigma_r^- = \left[%
\begin{array}{ccccc}
  0 & 0 & 1_r & 0 & 0 \\
  0 & 1_{n-1-r} & 0 & 0 & 0 \\
  1_r & 0 & 0 & 0 & 0 \\
  0 & 0 & 0 & 1_{n-1-r} & 0 \\
  0 & 0 & 0 & 0 & 1_2 \\
\end{array}%
\right] \in O^-(2n,q).
\]

For each $r$, with $0 \leq r \leq n-1$, put
\[
A_r^- = \{ w \in P^-(2n,q) \mid \sigma_r^- w(\sigma_r^-)^{-1} \in
P^-(2n,q)\}.
\]
As a disjoint union of right cosets of $P^- = P^-(2n,q)$, the
Bruhat decomposition in (30) can be written as
\begin{equation}
O^-(2n,q) = \coprod_{r=0}^{n-1} P^- \sigma_r^- (A_r^- \backslash
P^-).
\end{equation}

$Q^-(2n,q)$ is a subgroup of index 2 in $P^-(2n,q)$, defined by:
\begin{align*}
Q^- &= Q^-(2n,q) \\
    &= \left\{ \left[%
\begin{array}{ccc}
  A & 0 & 0\\
  0 & ^tA^{-1} & 0 \\
  0 & 0 & i \\
\end{array}%
\right] \left[%
\begin{array}{ccc}
  1_{n-1} & B & ^th^ti\eta i \\
  0 & 1_{n-1} & 0 \\
  0 & h & 1_2 \\
\end{array}%
\right] \Biggl \lvert %
\begin{array}{c}
  A \in GL(n-1,q), i \in SO^{-}(2,q),\\
  {}^tB + {}^th \delta_a h ~~\text{is alternating}\\
\end{array}%
\right\}.
\end{align*}
In fact, in view of (29), we have:
\[
P^-(2n,q) = Q^-(2n,q)\coprod \rho Q^-(2n,q),
\]with
\[
\rho = \left[%
\begin{array}{cccc}
  1_{n-1} & 0 & 0 & 0 \\
  0 & 1_{n-1} & 0 & 0 \\
  0 & 0 & 1 & 1 \\
  0 & 0 & 1 & 1 \\
\end{array}%
\right] \in P^-(2n,q).
\]

For each $r$, with $0 \leq r \leq n-1$, we define
\begin{align*}
B_r^- &= \{w \in Q^-(2n,q) \mid \sigma_r^- w (\sigma_r^-)^{-1} \in
P^-(2n,q)\}\\
&= \{w \in Q^-(2n,q) \mid \sigma_r^- w (\sigma_r^-)^{-1} \in
Q^-(2n,q)\},
\end{align*}
which is a subgroup of index 2 in $A_r^-$.

The decompositions in (30) and (31) can be modified so as to give:
\begin{equation}
O^-(2n,q) = \coprod_{r=0}^{n-1} P^- \sigma_r^- Q^- =
(\coprod_{r=0}^{n-1} Q^-\sigma_r^- Q^-)\coprod(\coprod_{r=0}^{n-1}
\rho Q^- \sigma_r^- Q^-),
\end{equation}
\begin{equation}
O^-(2n,q) = \coprod_{r=0}^{n-1} P^- \sigma_r^-(B_r^-\backslash
Q^-) = (\coprod_{r=0}^{n-1} Q^- \sigma_r^-(B_r^-\backslash Q^-))
\coprod (\coprod_{r=0}^{n-1} \rho Q^- \sigma_r^-(B_r^-\backslash
Q^-)).
\end{equation}
 The order of the general linear group  $GL(n,q)$ is given by
 \begin{equation}
g_n = \prod_{j=0}^{n-1} (q^n-q^j) = q^{{n \choose 2}} \prod
_{j=1}^n (q^j-1).
 \end{equation}
For integers $n,r$  with $0 \leq r \leq n$, the $q$-binomial
coefficients are defined as:
\[
\left[ \substack{n \\ r}
 \right]_q = \prod_{j=0}^{r-1} (q^{n-j} - 1)/(q^{r-j} - 1).
\]
Then, for integers $n,r$  with $0 \leq r \leq n$, we have
\begin{equation}
\frac{g_n}{g_{n-r} g_r} = q^{r(n-r)}\left[ \substack{n \\ r}
 \right]_q.
\end{equation}
In \cite{Kim8}, it is shown that
\begin{equation}
|A_r^-| = 2(q+1) g_r g_{n-1-r} q^{(n-1)(n+2)/2} q^{r(2n-3r-5)/2},
\end{equation}
\begin{equation}
|P^-(2n,q)| = 2(q+1) g_{n-1} q^{(n-1)(n+2)/2}.
\end{equation}
So, from (35)-(37), we get:
\begin{equation}
|A_r^- \backslash P^-(2n,q)| = |B_r^- \backslash Q^-(2n,q)| =
\left[ \substack{n-1 \\ r}
 \right]_q q^{r(r+3)/2},
\end{equation}
and
\begin{align}
\begin{split}
|Q^-(2n,q) \sigma_r^- Q^-(2n,q)| &= |\rho Q^-(2n,q) \sigma_r^-
Q^-(2n,q)|\\
&= \frac{1}{2}|P^-(2n,q) \sigma_r^- Q^-(2n,q)|\\
&= \frac{1}{2}|P^-(2n,q)||B_r^- \backslash Q^-(2n,q)|\\
&= \frac{1}{2}|P^-(2n,q)||A_r^- \backslash P^-(2n,q)|\\
&= \frac{1}{2}|P^-(2n,q)|^2|A_r^-|^{-1}\\
&= (q+1) q^{n^2-n} \prod_{j=1}^{n-1} (q^j-1)\left[ \substack{n-1\\
r} \right]_q q^{{r \choose 2}} q^{2r}
\end{split}
\end{align}(cf. (34), (37), (38)).

Let
\begin{align}
DC_1^+(n,q) &= Q^-(2n,q) \sigma_{n-1}^- Q^-(2n,q), \text{for $
n=2,4,6,\ldots$},\\
DC_2^+(n,q) &= Q^-(2n,q) \sigma_{n-2}^- Q^-(2n,q), \text{for $
n=2,4,6,\ldots$},\\
DC_3^+(n,q) &= \rho Q^-(2n,q) \sigma_{n-2}^- Q^-(2n,q), \text{for
$n=2,4,6,\ldots$},\\
DC_4^+(n,q) &= \rho Q^-(2n,q) \sigma_{n-3}^- Q^-(2n,q), \text{for
$n=4,6,8,\ldots$},\\
DC_1^-(n,q) &= Q^-(2n,q) \sigma_{n-1}^- Q^-(2n,q), \text{for $
n=1,3,5,\ldots$},\\
DC_2^-(n,q) &= Q^-(2n,q) \sigma_{n-2}^- Q^-(2n,q), \text{for $
n=3,5,7,\ldots$},\\
DC_3^-(n,q) &= \rho Q^-(2n,q) \sigma_{n-2}^- Q^-(2n,q), \text{for
$n=3,5,7,\ldots$},\\
DC_4^-(n,q) &= \rho Q^-(2n,q) \sigma_{n-3}^- Q^-(2n,q), \text{for
$n=3,5,7,\ldots$}.
\end{align}
Then, from (39), we have:
\begin{equation}
N_i^{\pm}(n,q) := |DC_i^{\pm}(n,q)| =
A_i^{\pm}(n,q)B_i^{\pm}(n,q), ~~\text{for $i=1,2,3,4$}
\end{equation}(cf. (1)-(16)).

\emph{Unless otherwise stated, from now on, we will agree that
anything related to $DC_1^+(n,q)$, $DC_2^+(n,q)$  and
$DC_3^+(n,q)$  are defined for $n=2,4,6,\ldots,$ anything related
to $DC_4^+(n,q)$ is defined for $n=4,6,8,\ldots,$ anything related
to $DC_1^-(n,q)$ is defined for $n=1,3,5,\ldots,$ and anything
related to $DC_2^-(n,q)$, $DC_3^-(n,q)$, and $DC_4^-(n,q)$ are
defined for $n=3,5,7,\ldots$.}

\section{Exponential sums over double cosets of $O^-(2n,2^r)$}

The following notations will be used throughout this paper:

\begin{gather*}
tr(x)=x+x^2+\cdots+x^{2^{r-1}} \text{the trace function} ~\fd_q
\rightarrow \fd_2,\\
\lambda(x) = (-1)^{tr(x)} ~\text{the canonical additive character
of} ~\fd_q.
\end{gather*}

Then any nontrivial additive character $\psi$ of $\fd_q$ is given
by $\psi(x) = \lambda(ax)$ , for a unique $a \in \fd_q^*$.

For any nontrivial additive character $\psi$ of $\fd_q$ and $a \in
\fd_q^*$, the Kloosterman sum  $K_{GL(t,q)}(\psi;a)$ for $GL(t,q)$
is defined as
\[
K_{GL(t,q)}(\psi;a) = \sum_{w \in GL(t,q)} \psi(Trw + aTrw^{-1}).
\]
Notice that, for $t=1$, $K_{GL(1,q)}(\psi;a)$ denotes the
Kloosterman sum $K(\psi;a)$.

For the Kloosterman sum $K(\psi ; a)$, we have the Weil bound (cf.
\cite{RH})
\begin{equation}
\mid K(\psi ; a) \mid \leq 2\sqrt{q}.
\end{equation}

 In \cite{Kim6}, it is shown that $K_{GL(t,q)}(\psi ; a)$ ~satisfies the following recursive relation:
  for integers $t \geq 2$, ~$a \in \fd_q^*$ ,
\begin{equation}
K_{GL(t,q)}(\psi ; a) = q^{t-1}K_{GL(t-1,q)}(\psi ; a)K(\psi ;a) +
q^{2t-2}(q^{t-1}-1)K_{GL(t-2,q)}(\psi ; a),
\end{equation}
where we understand that $K_{GL(0,q)}(\psi ; a)=1$ . From (51), in
\cite{Kim6} an explicit expression of the Kloosterman sum for
$GL(t,q)$ was derived.\\

\begin{theorem}[\cite{Kim6}] For integers
$t \geq 1$, and $a \in \fd_q^*$, the Kloosterman sum
$K_{GL(t,q)}(\psi ; a)$ is given by

\[
K_{GL(t,q)}(\psi ; a)=q^{(t-2)(t+1)/2} \sum_{l=1}^{[(t+2)/2]} q^l
K(\psi;a)^{t+2-2l}\sum \prod_{\nu=1}^{l-1} (q^{j_\nu -2\nu}-1),
\]
where  $K(\psi;a)$ is the Kloosterman sum and the inner sum is
over all integers $j_1,\ldots,j_{l-1}$ satisfying $2l-1 \leq
j_{l-1} \leq j_{l-2} \leq \cdots \leq j_1 \leq t+1$. Here we agree
that the inner sum is $1$ for $l=1$.
\end{theorem}\

\begin{proposition}[\cite{Kim7}, Prop. 3.1] Let  $\psi$ be a nontrivial additive character of $\fd_q$.
 Then
\begin{align}
 &(a) \sum_{i \in SO^-(2,q)} \psi(Tri) = K(\psi;1),~~~~~~~~~~~~~~~~~~~~~~~~~~~~~~~~~~~~~~~~~~~~ \\
 &(b) \sum_{i \in SO^-(2,q)} \psi(Tr\left[%
\begin{array}{cc}
  1 & 1 \\
  0 & 1 \\
\end{array}%
\right]i) = q+1.
 \end{align}
\end{proposition}

\begin{proposition}[\cite{Kim8}, Prop. 4.4] Let  $\psi$ be a nontrivial additive character of $\fd_q$.
 For each positive integer $r$, let $\Omega_r$ be the set of all  $r\times r$ nonsingular symmetric matrices over $\fd_q$.
 Then the $b_r(\psi)$ defined below is independent of $\psi$, and is equal to:
 \begin{align}
\begin{split}
b_r = b_r(\psi) &= \sum_{B \in \Omega_r} \sum_{h \in
\fd_q^{r\times 2}} \psi(Tr\delta_a{}^thBh)\\
&= \left\{%
\begin{array}{ll}
    q^{r(r+6)/4} \prod_{j=1}^{r/2} (q^{2j-1}-1), & \hbox{for $r$ even,} \\
    -q^{(r^2+4r-1)/4} \prod_{j=1}^{(r+1)/2} (q^{2j-1}-1), & \hbox{for $r$ odd.} \\
\end{array}%
\right.
\end{split}
 \end{align}
\end{proposition}

 In Section 5 of \cite{Kim8}, it is shown that the Gauss sum for $O^-(2n,q)$, with $\psi$ a nontrivial additive character of $\fd_q$,
   is given by:
\begin{align*}
\sum_{w \in O^-(2n,q)} \psi(Trw) &= \sum_{r=0}^{n-1} \sum_{w \in
P^-\sigma_r^-Q^-} \psi(Trw)\\
&= \sum_{r=0}^{n-1} \sum_{w \in Q^-\sigma_r^-Q^-} \psi(Trw) +
\sum_{r=0}^{n-1} \sum_{w \in \rho Q^-\sigma_r^-Q^-} \psi(Trw)
~~(cf. (32)),
\end{align*}
with
\begin{align}
\begin{split}
\sum_{w \in Q^-\sigma_r^-Q^-} \psi(Trw) = ~&|B_r^- \backslash
Q^-|\sum_{w \in Q^-} \psi(Trw\sigma_r^-)\\
=~& q^{(n-1)(n+2)/2} \sum_{i \in SO^-(2,q)} \psi(Tri)\\
&\times |B_r^- \backslash Q^-| q^{r(n-r-3)}b_r(\psi)
K_{GL(n-1-r,q)}(\psi;1),
\end{split}
\end{align}
\begin{align}
\begin{split}
\sum_{w \in \rho Q^-\sigma_r^-Q^-} \psi(Trw) = ~&|B_r^- \backslash
Q^-|\sum_{w \in Q^-} \psi(Tr\rho w\sigma_r^-)\\
=~& q^{(n-1)(n+2)/2} \sum_{i \in SO^-(2,q)} \psi(Tr\left[%
\begin{array}{cc}
  1 & 1 \\
  0 & 1 \\
\end{array}%
\right]i)\\
&\times |B_r^- \backslash Q^-| q^{r(n-r-3)}b_r(\psi)
K_{GL(n-1-r,q)}(\psi;1).
\end{split}
\end{align}
Here one uses (33) and the fact that  $\rho^{-1} w \rho \in Q^-$,
for all $w \in Q^-$.

 Now, we see from (52)-(56) and (38) that, for each $r$  with $0
 \leq r \leq n-1$,
\begin{align}
\begin{split}
\sum_{w \in Q^-\sigma_r^-Q^-} \psi(Trw) = ~&q^{(n-1)(n+2)/2}\left[
\substack{n-1 \\ r} \right]_q K(\psi;1) K_{GL(n-1-r,q)}(\psi;1) \\
& \times \left\{%
\begin{array}{ll}
    -q^{rn-\frac{1}{4}r^2}\prod_{j=1}^{r/2} (q^{2j-1}-1), & \hbox{for $r$ even,} \\
    q^{rn-\frac{1}{4}(r+1)^2} \prod_{j=1}^{(r+1)/2} (q^{2j-1}-1), & \hbox{for $r$ odd,} \\
\end{array}%
\right.
\end{split}
\end{align}

\begin{align}
\begin{split}
\sum_{w \in \rho Q^-\sigma_r^-Q^-} \psi(Trw) =
~&(q+1)q^{(n-1)(n+2)/2}\left[
\substack{n-1 \\ r} \right]_q K_{GL(n-1-r,q)}(\psi;1) \\
& \times \left\{%
\begin{array}{ll}
    q^{rn-\frac{1}{4}r^2}\prod_{j=1}^{r/2} (q^{2j-1}-1), & \hbox{for $r$ even,} \\
    -q^{rn-\frac{1}{4}(r+1)^2} \prod_{j=1}^{(r+1)/2} (q^{2j-1}-1), & \hbox{for $r$ odd.} \\
\end{array}%
\right.
\end{split}
\end{align}

For our purposes, we need the following special cases of
exponential sums in (57) and (58). We state them separately as a
theorem.

\begin{theorem}
Let $\psi$  be any nontrivial additive character of $\fd_q$.
Then, in the notations of  (1), (3), (5), (7), (9), (11), (13),
and (15), we have
\begin{align*}
\sum_{w \in DC_i^\pm(n,q)} \psi(Trw) &= \pm A_i^\pm(n,q)
K(\psi;1), \text{for $i=1,3$},\\
\sum_{w \in DC_2^\pm(n,q)} \psi(Trw) &= \pm(-1) A_2^\pm(n,q)
K(\psi;1)^2,\\
\sum_{w \in DC_4^\pm(n,q)} \psi(Trw) &= \pm(-1) q^{-1}
A_4^\pm(n,q) K_{GL(2,q)}(\psi;1)\\
&= \pm(-1) A_4^\pm(n,q) (K(\psi;1)^2 + q^2 - q)
\end{align*}(cf. (40)-(47), (51)).
\end{theorem}

\begin{proposition}[\cite{Kim2}] For $n=2^s(s \in \Z_{\geq 0})$, and  $\psi$ a nontrivial additive character of $\fd_q$,
\[
K(\psi;a^n) = K(\psi;a).
\]
\end{proposition}

 For the next corollary, we need a result of Carlitz.

 \begin{theorem}[\cite{L2}]  For the canonical additive character $\lambda$ of $\fd_q$, and $a \in \fd_q^*$,
 \begin{equation}
K_2(\lambda;a) = K(\lambda;a)^2 - q.
 \end{equation}
 \end{theorem}

The next corollary follows from Theorems 6 and 8, Proposition 7,
and by simple change of variables.

\begin{corollary}
Let $\lambda$  be the canonical additive character of $\fd_q$, and
let $a \in \fd_q^*$. Then we have
\begin{align}
\sum_{w \in DC_i^\pm(n,q)} \lambda(aTrw) &= \pm A_i^\pm(n,q)
K(\lambda;a), \text{for $i=1,3$},\\
\begin{split}
\sum_{w \in DC_2^\pm(n,q)} \lambda(aTrw) &= \pm(-1) A_2^\pm(n,q)
K(\lambda;a)^2\\
&= \pm(-1) A_2^\pm(n,q) (K_2(\lambda;a) + q),
\end{split}
\end{align}
\begin{align}
\begin{split}
\sum_{w \in DC_4^\pm(n,q)} \lambda(aTrw) &= \pm(-1) A_4^\pm(n,q)
(K(\lambda;a)^2 + q^2 - q)\\
&= \pm(-1) A_4^\pm(n,q) (K_2(\lambda;a) + q^2).
\end{split}
\end{align}
\end{corollary}

\begin{proposition}[\cite{Kim2}]
 Let  $\lambda$ be the canonical additive character of $\fd_q$, $m \in \Z_{>0}$, $\beta \in \fd_q$. Then
\begin{align}
\begin{split}
& \sum_{a \in \fd_q^*} \lambda(-a \beta) K_m(\lambda;a)\\
&= \left\{%
\begin{array}{ll}
    q K_{m-1}(\lambda ; \beta^{-1}) + (-1)^{m+1}, & \hbox{if $\beta \neq 0$,} \\
    (-1)^{m+1}, & \hbox{if $\beta = 0$,} \\
\end{array}%
\right.
\end{split}
\end{align}with the convention $K_0(\lambda;\beta^{-1}) =
\lambda(\beta^{-1}).$
\end{proposition}

 For any integer $r$  with $0 \leq r \leq n-1$, and each $\beta \in \fd_q$,  we let
 \begin{align*}
N_{Q^-\sigma_r^-Q^-}(\beta) &= \mid \{w \in Q^-\sigma_r^-Q^- \mid
Trw = \beta \}\mid,\\
N_{\rho Q^-\sigma_r^-Q^-}(\beta) &= \mid \{w \in \rho
Q^-\sigma_r^-Q^- \mid Trw = \beta \}\mid.
 \end{align*}
Then it is easy to see that
\begin{align}
q N_{Q^-\sigma_r^-Q^-}(\beta) &= |Q^-\sigma_r^-Q^-| + \sum_{a \in
\fd_q^*} \lambda(-a\beta) \sum_{w \in Q^-\sigma_r^-Q^-}
\lambda(aTrw),\\
q N_{\rho Q^-\sigma_r^-Q^-}(\beta) &= |\rho Q^-\sigma_r^-Q^-| +
\sum_{a \in \fd_q^*} \lambda(-a\beta) \sum_{w \in \rho
Q^-\sigma_r^-Q^-} \lambda(aTrw).
\end{align}

Now, from (60)-(65) and (40)-(48), we have the following result.

\begin{proposition}
(a)  For $i=1,3,$
\begin{equation}
N_{DC_i^\pm(n,q)}(\beta) = q^{-1} A_i^\pm(n,q) B_i^\pm(n,q) \pm
q^{-1} A_i^\pm(n,q) \times \left\{%
\begin{array}{ll}
    1, & \hbox{$\beta = 0$,} \\
    q+1, & \hbox{$tr(\beta^{-1}) = 0$,} \\
    -q+1, & \hbox{$tr(\beta^{-1}) = 1$,} \\
\end{array}%
\right.
\end{equation}
\begin{multline}
(b)~~ N_{DC_2^\pm(n,q)}(\beta) = q^{-1} A_2^\pm(n,q) B_2^\pm(n,q)
\pm (-1) q^{-1} A_2^\pm(n,q) \\
\times \left\{%
\begin{array}{ll}
    qK(\lambda;\beta^{-1})-q-1, & \hbox{$\beta \neq 0$,} \\
    q^2-q-1, & \hbox{$\beta = 0$,} \\
\end{array}%
\right.
\end{multline}
\begin{multline}
(c) ~~N_{DC_4^\pm(n,q)}(\beta) = q^{-1} A_4^\pm(n,q) B_4^\pm(n,q)
\pm (-1) q^{-1} A_4^\pm(n,q) \\
\times \left\{%
\begin{array}{ll}
    qK(\lambda;\beta^{-1})-q^2-1, & \hbox{$\beta \neq 0$,} \\
    q^3-q^2-1, & \hbox{$\beta = 0$.} \\
\end{array}%
\right.
\end{multline}
\end{proposition}

\begin{corollary}
(a) For all even $n \geq 2$ and all $q$, $N_{DC_i^+(n,q)}(\beta) >
0$, for all $\beta$ and $i=1,2.$\\
(b) For all even $n \geq 4$  and all $q$, $N_{DC_3^+(n,q)}(\beta)
> 0$, for all $\beta$; for $n=2$ and all $q$,
\begin{equation}
N_{DC_3^+(2,q)}(\beta) = \left\{%
\begin{array}{ll}
    q^3+q^2, & \hbox{$\beta = 0$,} \\
    2q^3 + 2q^2, & \hbox{$tr(\beta^{-1}) = 0$,} \\
    0, & \hbox{$tr(\beta^{-1})=1$.} \\
\end{array}%
\right.
\end{equation}
(c) For all even $n \geq 4$ and all $q$,
$N_{DC_4^+(n,q)}(\beta)>0$, for all $\beta$.\\
(d) For all odd $n \geq 3$  and all  $q$,
$N_{DC_1^-(n,q)}(\beta)>0$, for all $\beta$;  for $n=1$ and all
$q$,
\begin{equation}
N_{DC_1^-(1,q)}(\beta) = \left\{%
\begin{array}{ll}
    1, & \hbox{$\beta = 0$,} \\
    0, & \hbox{$tr(\beta^{-1}) = 0$,} \\
    2, & \hbox{$tr(\beta^{-1}) = 1$.} \\
\end{array}%
\right.
\end{equation}
(e)  For all odd $n \geq 3$ and all $q$,
$N_{DC_i^-(n,q)}(\beta)>0$, for all $\beta$  and $i=2,3.$\\
(f) For all odd $n \geq 5$ and all $q$, or $n=3$ and all $q \geq
4$, $N_{DC_4^-(n,q)}(\beta)>0$, for all $\beta$; for $n=3$ and
$q=2$,
\begin{equation}
N_{DC_4^-(3,2)}(\beta) = \left\{%
\begin{array}{ll}
    576 = |\rho Q^-(6,2)|, & \hbox{$\beta = 0$,} \\
    0, & \hbox{$\beta = 1$.} \\
\end{array}%
\right.
\end{equation}
\end{corollary}
\begin{pf*}{Proof.} All assertions except (f) are left to the
reader.\\
(f) Let $\beta = 0$. Then $N_{DC_4^-(n,q)}(0)>0$, for all odd $n
\geq 3$ and all $q$, as one can see from (68). Now, let $\beta
\neq 0$.  Then, by invoking the Weil bound in (50), we have
\begin{multline}
N_{DC_4^-(n,q)}(\beta) \geq q^{-1} A_4^-(n,q)\\
\times \{ q^{\frac{1}{4}(n-3)^2}(q^{n-2}-1) \prod_{j=1}^{(n-1)/2}
(q^{2j}-1)-(q^2+2q^{\frac{3}{2}}+1) \}.
\end{multline}
Let $n \geq 5$.  Then we see from (72) that, for all $q$,
\[
N_{DC_4^-(n,q)}(\beta) \geq q^{-1}
A_4^-(n,q)\{q(q^3-1)-(q^2+2q^{\frac{3}{2}} + 1)\}>0.
\]
If  $n=3$ and $q \geq 4$, then, from (72), we have
\[
N_{DC_4^-(3,q)}(\beta) \geq q^{-1}
A_4^-(3,q)\{(q-1)(q^2-1)-(q^2+2q^{\frac{3}{2}} + 1)\}>0.
\]
On the other hand, if  $n=3$ and $q=2$, then we get the values in
(71) directly from
(68).~~~~~~~~~~~~~~~~~~~~~~~~~~~~~~~~~~~~~~~~~~~~~~~~~~~~~~~~~~~~~~~~~~~~~~~~~~~~~~~~~~~~~~~~~~~~~~~~~~~~~~~~~~~~~~~~~~~~~$\square$
\end{pf*}

\section{Construction of codes}

 Here we will construct eight infinite families of binary linear codes $C(DC_1^+(n,q))$
 of length $N_1^+(n,q)$, $C(DC_2^+(n,q))$  of length $N_2^+(n,q)$, $C(DC_3^+(n,q))$
  of length $N_3^+(n,q)$, for $n=2,4,6,\ldots$ and all $q$; $C(DC_4^+(n,q))$  of length
  $N_4^+(n,q)$, for $n=4,6,8,\ldots$ and all $q$; $C(DC_1^-(n,q))$  of length $N_1^-(n,q)$ for $n=1,3,5,\ldots$
   and all $q$; $C(DC_2^-(n,q))$ of length $N_2^-(n,q)$, $C(DC_3^-(n,q))$  of length  $N_3^-(n,q)$, $C(DC_4^-(n,q))$
   of length $N_4^-(n,q)$, for $n=3,5,7,\ldots$ and all $q$, respectively associated with the double cosets
   $DC_1^+(n,q)$, $DC_2^+(n,q)$, $DC_3^+(n,q)$, $DC_4^+(n,q)$, $DC_1^-(n,q)$, $DC_2^-(n,q)$, $DC_3^-(n,q)$, $DC_4^-(n,q)$ (cf. (40)-(48)).

 Let $g_1,g_2,\ldots,g_{N_i^\pm(n,q)}$ be  fixed orderings of the elements in  $DC_i^\pm(n,q)$, for $i=1,2,3,4,$
 by abuse of notations.  Then we put
 \[
v_i^\pm(n,q) = (Trg_1, Trg_2,\cdots, Trg_{N_i^\pm(n,q)}) \in
\fd_q^{N_i^\pm(n,q)}, \text{for $i=1,2$}.
 \]

 The binary codes $C(DC_1^+(n,q))$, $C(DC_2^+(n,q))$,
 $C(DC_3^+(n,q))$, $C(DC_4^+(n,q))$, $C(DC_1^-(n,q))$,
 $C(DC_2^-(n,q))$, $C(DC_3^-(n,q))$, and $C(DC_4^-(n,q))$ are
 defined as:
 \begin{equation}
C(DC_i^\pm(n,q)) = \{u \in \fd_q^{N_i^\pm(n,q)} \mid u\cdot
v_i^\pm(n,q) = 0 \}, \text{for}~~ i=1,2,3,4,
 \end{equation}
where the dot denotes respectively the usual inner product in
$\fd_q^{N_i^\pm(n,q)}$, for $i=1,2,3,4.$

The following theorem of  Delsarte  is well-known.
\begin{theorem}[\cite{FJ}]
Let  $B$ be a linear code over $\fd_q$.  Then
\[
(B|_{\fd_2})^\bot = tr(B^\bot).
\]
\end{theorem}

  In view of this theorem, the respective duals of the codes in (73) are given by:
  \begin{multline}
C(DC_i^\pm(n,q))^\bot=\{c_i^\pm(a)=c_i^\pm(a;n,q)=(tr(aTrg_1),
\ldots, tr(aTrg_{N_i^\pm(n,q)}) \mid a \in \fd_q\},
\end{multline} for $i=1,2,3,4.$

  Let  $\fd_2^+, \fd_q^+$ denote the additive groups of the fields $\fd_2,\fd_q$, respectively.
  Then we have the following exact sequence of groups:
\begin{equation*}
0 \rightarrow \fd_2^+ \rightarrow \fd_q^+ \rightarrow
\Theta(\fd_q) \rightarrow 0,
\end{equation*}
where the first map is the inclusion and the second one is  the
Artin-Schreier operator in characteristic two given by $\Theta(x)
= x^2+x$.   So
\begin{equation}
\Theta(\fd_q) = \{\alpha^2 + \alpha \mid  \alpha \in \fd_q \},~
and ~~[\fd_q^+ : \Theta(\fd_q)] = 2.
\end{equation}

\begin{theorem}[\cite{Kim2}]
 Let $\lambda$  be the canonical additive character of $\fd_q$, and let $\beta \in \fd_q^*$. Then
\begin{align}
 &(a) \sum_{\alpha \in
 \fd_q-\{0,1\}}\lambda(\frac{\beta}{\alpha^2+\alpha})= K(\lambda;\beta)-1,~~~~~~~~~~~~~~~~~~~~~~~~~~~~~~~~~~~~~~~\\
 &(b)\sum_{\alpha \in
\fd_q}\lambda(\frac{\beta}{\alpha^2+\alpha+b}) =
-K(\lambda;\beta)-1,
\end{align}
if $x^2+x+b (b \in \fd_q)$ is irreducible over $\fd_q$, or
equivalently if $b \in \fd_q\setminus\Theta(\fd_q)$ (cf.(75)).
\end{theorem}

\begin{theorem}
(a) The map $\fd_q \rightarrow C(DC_i^+(n,q))^\bot$ ($a \mapsto
c_i^+(a)$) ($i=1,2$) is an $\fd_2$-linear isomorphism for $n
\geq 2$ even and all $q$.\\
(b) The map $\fd_q \rightarrow C(DC_3^+(n,q))^\bot$ ($a \mapsto
c_3^+(a)$) is an $\fd_2$-linear isomorphism for $n
\geq 4$ even and all $q$, or $n=2$ and $q \geq 8$.\\
(c) The map $\fd_q \rightarrow C(DC_4^+(n,q))^\bot$ ($a \mapsto
c_4^+(a)$) is an $\fd_2$-linear isomorphism for $n \geq 4$ even
and all $q$.\\
(d) The map $\fd_q \rightarrow C(DC_1^-(n,q))^\bot$ ($a \mapsto
c_1^-(a)$) is an $\fd_2$-linear isomorphism for $n \geq 1$ odd and
all $q$.\\
(e) The map $\fd_q \rightarrow C(DC_i^-(n,q))^\bot$ ($a \mapsto
c_i^-(a)$)($i=2,3$) is an $\fd_2$-linear isomorphism for $n \geq
3$ odd and all $q$.\\
(f) The map $\fd_q \rightarrow C(DC_4^-(n,q))^\bot$ ($a \mapsto
c_4^-(a)$) is an $\fd_2$-linear isomorphism for $n \geq 5$ odd and
all $q$, or $n=3$ and $q \geq 4$.\\
\end{theorem}
\begin{pf*}{Proof.}
All maps are clearly $\fd_2$-linear and surjective. Let $a$ be in
the kernel of map $\fd_q \rightarrow C(DC_1^+(n,q))^\bot$ ($a
\mapsto c_1^+(a)$). Then $tr(aTrg)=0$, for all $g \in
DC_1^+(n,q)$.  Since, by Corollary 12(a), $Tr : DC_1^+(n,q)
\rightarrow \fd_q$ is surjective, $tr(a\alpha)=0$, for all $\alpha
\in \fd_q$. This implies that $a=0$, since otherwise $tr : \fd_q
\rightarrow \fd_2$ would be the zero map. This shows (a). All the
other assertions can be handled in the same way, except for $n=2$
and $q \geq 8$ case  of (b) and $n=1$ case of (d). Assume first
that we are in  the $n=2$ and $q \geq 8$ case  of  (b). Let $a$ be
in the kernel of the map $\fd_q \rightarrow C(DC_3^+(2,q))^\bot$
($a \mapsto c_3^+(a)$). Then, by (69), $tr(a\beta)=0$, for all
$\beta \in \fd_q^*$, with $tr(\beta^{-1})=0$.  Hilbert's theorem
90 says that $tr(\gamma)=0$ $\Leftrightarrow$ $\gamma =
\alpha^2+\alpha$, for some $\alpha \in \fd_q$, and hence
$\sum_{\alpha \in
 \fd_q-\{0,1\}}\lambda(\frac{a}{\alpha^2+\alpha}) = q-2$. If $a
 \neq 0$, then, using (76) and the Weil bound (50), we would have
 \[
 q-2 = \sum_{\alpha \in
 \fd_q-\{0,1\}}\lambda(\frac{a}{\alpha^2+\alpha}) = K(\lambda;a) -
 1 \leq 2\sqrt{q} - 1.
 \]But this is impossible, since  $x > 2\sqrt{x}+1$, for $x \geq
 8$.\\
Assume next that we are in the $n=1$ case of  (d).   Let $a$ be in
the kernel of the map $\fd_q \rightarrow C(DC_1^-(1,q))^\bot$($a
\mapsto c_1^-(a)$). Then, by (70), $tr(a\beta)=0$, for all $\beta
\in \fd_q^*$  with $tr(\beta^{-1})=1$.  Let $b \in \fd_q
\backslash \Theta(\fd_q)$. Then $tr(\gamma)=1$ $\Leftrightarrow$
$\gamma = \alpha^2+\alpha+b$, for some $\alpha \in \fd_q$. As
$z^2+z+b$ is irreducible over $\fd_q$, $\alpha^2+\alpha+b \neq 0$,
for all $\alpha \in \fd_q$, and hence
$tr(\frac{a}{\alpha^2+\alpha+b})=0$, for all $\alpha \in \fd_q$.
So $\dis\sum_{\alpha \in
\fd_q}\lambda(\frac{a}{\alpha^2+\alpha+b}) = q$. Assume now that
$a \neq 0$. Then, from (77) and (50),
\[
q = -K(\lambda;a)-1 \leq 2\sqrt{q}-1.
\]But this is impossible, since , $x > 2\sqrt{x}-1$, for $x \geq
2$.~~~~~~~~~~~~~~~~~~~~~~~~~~~~~~$\square$
\end{pf*}

Remark : One can show that the kernel of the maps $\fd_q
\rightarrow C(DC_3^+(2,q))^\bot$ ($a \mapsto c_3^+(a)$), for
$q=2,4$, and of the map $\fd_2 \rightarrow C(DC_4^-(3,2))^\bot$
($a \mapsto c_4^-(a)$) are all equal to $\fd_2$.

\section{Recursive formulas for power moments of Kloosterman sums}

Here we will be able to find, via Pless power moment identity,
infinite families of recursive formulas generating power moments
of Kloosterman and 2-dimensional Kloosterman sums over all $\fd_q$
(with three exceptions) in terms of the frequencies of weights in
$C(DC_i^\pm(n,q))$, for $i=1,3$ and $C(DC_i^\pm(n,q))$, for
$i=2,4$, respectively.

\begin{theorem}[Pless power moment identity, \cite{FJ}]
Let $B$ be an $q$-ary $[n,k]$ code, and let $B_i$ (resp.
$B_i^\bot$) denote the number of codewords of weight $i$ in
$B$(resp. in $B^\bot$).  Then, for $h=0,1,2,\ldots,$
\begin{equation}
\sum_{j=0}^n j^h B_j=\sum_{j=0}^{min\{n,h\}} (-1)^j B_j^\bot
 \sum_{t=j}^h t! S(h,t) q^{k-t}(q-1)^{t-j}{n-j\choose n-t},
\end{equation}
where $S(h,t)$ is the Stirling number of the second kind defined
in (19).
\end{theorem}

\begin{lemma}
Let $c_i^\pm(a) = (tr(Trg_1),\ldots,tr(Trg_{N_i^\pm(n,q)})) \in
C(DC_i^\pm(n,q))^\bot$, for $a \in \fd_q^*$ and $i=1,2,3,4.$ Then
their Hamming weights are expressed as follows:
\begin{align}
(a) w(c_i^\pm(a)) &=  \frac{1}{2} A_i^\pm(n,q) \{B_i^\pm(n,q) \pm
(-1) K(\lambda;a)\}, ~\text{for}~~ i=1,3,~~~~~~\\
(b) w(c_2^\pm(a)) &= \frac{1}{2} A_2^\pm(n,q) (B_2^\pm(n,q) \pm
K(\lambda;a)^2)\\
&= \frac{1}{2} A_2^\pm(n,q) \{B_2^\pm(n,q) \pm (q +
K_2(\lambda;a))
\},\\
(c) w(c_4^\pm(a)) &= \frac{1}{2} A_4^\pm(n,q) \{B_4^\pm(n,q) \pm
(q^2 - q + K(\lambda;a)^2)\} \\
&=  \frac{1}{2} A_4^\pm(n,q) \{B_4^\pm(n,q) \pm (q^2 +
K_2(\lambda;a))\}
\end{align}(cf. (1)-(16)).
\end{lemma}
\begin{pf*}{Proof}
$ w(c_i^\pm(a)) =  \dis\frac{1}{2} \sum_{j=1}^{N_i^\pm(n,q)}
(1-(-1)^{tr(aTrg_j)}) = \frac{1}{2}(N_i^\pm(n,q) - \dis\sum_{w \in
DC_i^\pm(n,q)} \lambda(aTrw))$, for $i=1,2,3,4$. Our results now
follow from (48) and
(59)-(62).~~~~~~~~~~~~~~~~~~~~~~~~~~~~~~~$\square$
\end{pf*}

Let $u = (u_1, \ldots, u_{N_i^\pm(n,q)}) \in
\fd_2^{N_i^\pm(n,q)}$, for $i=1,2,3,4,$ with $\nu_\beta$ 1's in
the coordinate places where $Tr(g_j)=\beta$, for each $\beta \in
\fd_q$. Then from the definition of the codes $C(DC_i^\pm(n,q))$
(cf. (73)) we see that $u$ is a codeword with weight $j$ if and
only if $\dis\sum_{\beta \in \fd_q} \nu_\beta = j$ and $\dis
\sum_{\beta \in \fd_q} \nu_\beta \beta = 0$ (an identity in
$\fd_q$).  As there are $\dis \prod_{\beta \in \fd_q}
{N_{DC_i^\pm(n,q)}(\beta) \choose \nu_\beta}$ many such codewords
with weight $j$, we obtain the following result.

\begin{proposition}
Let $\{C_{i,j}^\pm(n,q)\}_{j=0}^{N_i^\pm(n,q)}$ be the weight
distribution of $C(DC_i^\pm(n,q))$, for $i=1,2,3,4$. Then we have
\begin{equation}
C_{i,j}^\pm(n,q) = \sum \prod_{\beta \in \fd_q}
{N_{DC_i^\pm(n,q)}(\beta) \choose \nu_\beta}, ~\text{for} ~~0 \leq
j \leq N_i^\pm(n,q), \text{and} ~~i=1,2,3,4,
\end{equation}
where the sum is over all the sets of integers
$\{\nu_\beta\}_{\beta \in \fd_q}$ $(0 \leq \nu_\beta \leq
N_{DC_i^\pm(n,q)}(\beta))$, satisfying
\begin{equation}
\sum_{\beta \in \fd_q} \nu_\beta = j, ~~\text{and}~~ \sum_{\beta
\in \fd_q} \nu_\beta \beta = 0.
\end{equation}
\end{proposition}

\begin{corollary}
Let $\{C_{i,j}^\pm(n,q)\}_{j=0}^{N_i^\pm(n,q)}$ be the weight
distribution of $C(DC_i^\pm(n,q))$, for $i=1,2,3,4$. Then we have
\[
C_{i,j}^\pm(n,q) = C_{i,N_i^\pm(n,q)-j}^\pm(n,q), ~~\text{for all
~~$j$, ~~with ~~$0 \leq j \leq N_i^\pm(n,q)$}.
\]
\end{corollary}
\begin{pf*}{Proof}
Under the replacements $\nu_\beta \rightarrow
N_{DC_i^\pm(n,q)}(\beta)-\nu_\beta$, for each $\beta \in \fd_q$,
the first equation in (85) is changed to $N_i^\pm(n,q)-j$, while
the second one in there and the summands in (84) are left
unchanged. The second sum in (85) is left unchanged, since $\dis
\sum_{\beta \in \fd_q} N_{DC_i^\pm(n,q)}(\beta)\beta = 0$, as one
can see by using the explicit expressions of
$N_{DC_i^\pm(n,q)}(\beta)$ in
(66)-(68).~~~~~~~~~~~~~~~~~~~~~~~~~~~~~~~~~~~~~~~~~~~~~~$\square$
\end{pf*}

\begin{theorem}[\cite{LJ}]
Let $q=2^r$, with $r \geq 2$. Then the range $R$ of
$K(\lambda;a)$, as a varies over $\fd_q^*$, is given by:
\[
R = \{ \tau \in \Z \mid |\tau|< 2\sqrt{q}, \tau \equiv -1 (mod
~4)\}.
\]
In addition, each value $\tau \in R$ is attained exactly
$H(\tau^2-q)$ times, where $H(d)$ is the Kronecker class number of
$d$.
\end{theorem}

The formulas appearing in the next theorem and stated in (18),
(22), and (25) follow  by applying the formula in (84) to each
$C(DC_i^\pm(n,q))$, using the explicit values of
$N_{DC_i^\pm(n,q)}(\beta)$ in (66)-(68), and taking Theorem 20
into consideration.

\begin{theorem}
Let $\{C_{i,j}^\pm(n,q)\}_{j=0}^{N_i^\pm(n,q)}$ be the weight
distribution of $C(DC_i^\pm(n,q))$, for $i=1,2,3,4$, and assume
that $q \geq 4$, for $C(DC_i^\pm(n,q))$($i=2,4$). Then we have\\
(a) For $i=1,3$, and $j=0,\ldots,N_i^\pm(n,q)$,
\begin{multline*}
C_{i,j}^\pm(n,q) = \sum {q^{-1} A_i^\pm(n,q)(B_i^\pm(n,q) \pm 1)
\choose \nu_0} \\
\times \prod_{tr(\beta^{-1})=0} {q^{-1} A_i^\pm(n,q)(B_i^\pm(n,q)
\pm (q+1)) \choose \nu_\beta} \prod_{tr(\beta^{-1})=1} {q^{-1}
A_i^\pm(n,q)(B_i^\pm(n,q) \pm (-q+1)) \choose \nu_\beta},
\end{multline*}
where the sum is over all the sets of nonnegative integers
$\{\nu_\beta\}_{\beta \in \fd_q}$ satisfying  $\sum_{\beta \in
\fd_q} \nu_\beta = j$, and $\sum_{\beta \in \fd_q} \nu_\beta \beta
= 0$.\\

(b) For $j=0,\ldots,N_2^\pm(n,q)$,
\begin{multline*}
C_{2,j}^\pm(n,q) = \sum {q^{-1} A_2^\pm(n,q)(B_2^\pm(n,q) \pm
(q+1-q^2)) \choose \nu_0}\\
\times \prod_{\substack{|\tau|<2\sqrt{q}\\ \tau \equiv-1(mod 4)}}
\prod_{K(\lambda;\beta^{-1})=\tau} {q^{-1}
A_2^\pm(n,q)(B_2^\pm(n,q) \pm (q+1-q\tau)) \choose \nu_\beta},
\end{multline*}
where the sum is over all the sets of nonnegative integers
$\{\nu_\beta\}_{\beta \in \fd_q}$ satisfying  $\sum_{\beta \in
\fd_q} \nu_\beta = j$, and $\sum_{\beta \in \fd_q} \nu_\beta \beta
= 0$.\\

(c) For $j=0,\ldots,N_4^\pm(n,q)$,
\begin{multline*}
 C_{4,j}^\pm(n,q) = \sum {q^{-1}
A_4^\pm(n,q) (B_4^\pm(n,q) \pm
(q^2+1-q^3)) \choose  \nu_0} \\
\times \prod_{\substack{|\tau|<2\sqrt{q}\\ \tau \equiv-1(mod 4)}}
\prod_{K(\lambda;\beta^{-1})=\tau} {q^{-1}
A_4^\pm(n,q)(B_4^\pm(n,q) \pm (q^2+1-q\tau)) \choose \nu_\beta},
\end{multline*}
where the sum is over all the sets of nonnegative integers
$\{\nu_\beta\}_{\beta \in \fd_q}$ satisfying  $\sum_{\beta \in
\fd_q} \nu_\beta = j$, and $\sum_{\beta \in \fd_q} \nu_\beta \beta
= 0$.
\end{theorem}

From now on, we will assume that, for  $C(DC_1^+(n,q))^\bot$,
$n\geq 2$ even and all $q$;  for $C(DC_2^+(n,q))^\bot$, $n \geq 2$
even and $q \geq 4$; for $C(DC_3^+(n,q))^\bot$, either $n \geq 4$
even and all $q$, or $n=2, q \geq 8$; for $C(DC_4^+(n,q))^\bot$,
$n \geq 4$ even and $q \geq 4$; for $C(DC_1^-(n,q))^\bot$, $n \geq
1$ odd and all $q$;  for $C(DC_2^-(n,q))^\bot$, $n \geq 3$ odd and
$q \geq 4$; for $C(DC_3^-(n,q))^\bot$, $n \geq 3$ odd and all $q$;
for $C(DC_4^-(n,q))^\bot$, $n \geq 3$ odd and $q \geq 4$. Under
these assumptions, each codeword in $C(DC_i^\pm(n,q))^\bot$ can be
written as $c_i^\pm(a)$, for $i=1,2,3,4$, and a unique $a \in
\fd_q$ (cf. Theorem 15, (74)).

Now, we apply the Pless power moment identity in (78) to
$C(DC_i^\pm(n,q))^\bot$, for those values of $n$ and $q$, in order
to get the results in Theorem 1(cf. (17), (18), (20)-(25)) about
recursive formulas.

The left hand side of that identity in (78) is equal to
\[
\sum_{a \in \fd_q^*} w(c_i^\pm(a))^h,
\]with  $w(c_i^\pm(a))$ given by (79)-(83). We have, for $i=1,3$,
\begin{align}
\begin{split}
\sum_{a \in \fd_q^*} w(c_i^\pm(a))^h &=
\frac{1}{2^h}A_i^\pm(n,q)^h \sum_{a \in \fd_q^*} \{B_i^\pm(n,q)
\pm (-1) K(\lambda;a)\}^h\\
&= \frac{1}{2^h}A_i^\pm(n,q)^h \sum_{l=0}^h (\pm(-1))^l {h \choose
l} B_i^\pm(n,q)^{h-l} MK^l.
\end{split}
\end{align}
Similarly, we have
\begin{align}
\sum_{a \in \fd_q^*} w(c_2^\pm(a))^h &=
\frac{1}{2^h}A_2^\pm(n,q)^h \sum_{l=0}^h (\pm1)^l {h \choose l}
B_2^\pm(n,q)^{h-l} MK^{2l} \\
&= \frac{1}{2^h}A_2^\pm(n,q)^h \sum_{l=0}^h (\pm1)^l {h \choose l}
(B_2^\pm(n,q) \pm q)^{h-l} MK_2^l,\\
\sum_{a \in \fd_q^*} w(c_4^\pm(a))^h &=
\frac{1}{2^h}A_4^\pm(n,q)^h \sum_{l=0}^h (\pm1)^l {h \choose l}
\{B_4^\pm(n,q) \pm (q^2-q)\}^{h-l} MK^{2l}\\
&= \frac{1}{2^h}A_4^\pm(n,q)^h \sum_{l=0}^h (\pm1)^l {h \choose l}
(B_4^\pm(n,q) \pm q^2)^{h-l} MK_2^l.
\end{align}

Note here that, in view of (59), obtaining power moments of
2-dimensional Kloosterman sums is equivalent to getting even power
moments of Kloosterman sums. Also, one has to separate the term
corresponding to $l=h$ in (86)-(90), and notes
$dim_{\fd_2}C(DC_i^\pm(n,q))^\bot = r$.

\begin{flushleft}
{\textbf{Acknowledgments}}
\end{flushleft}
This work was supported by grant No. R01-2008-000-11176-0 from the
Basic Research Program of the Korea Science and Engineering
Foundation.

\end{document}